\documentclass[12pt]{article}
\usepackage{amscd2000, amsmath2000}
\begin{document} 
\renewcommand{\thesubsection}{\arabic{subsection}}
\newenvironment{eq}{\begin{equation}}{\end{equation}}
\newenvironment{proof}{{\bf Proof}:}{\vskip 5mm }
\newenvironment{rem}{{\bf Remark}:}{\vskip 5mm }
\newenvironment{remarks}{{\bf Remarks}:\begin{enumerate}}{\end{enumerate}}
\newenvironment{examples}{{\bf Examples}:\begin{enumerate}}{\end{enumerate}}  
\newtheorem{proposition}{Proposition}[subsection]
\newtheorem{lemma}[proposition]{Lemma}
\newtheorem{definition}[proposition]{Definition}
\newtheorem{theorem}[proposition]{Theorem}
\newtheorem{cor}[proposition]{Corollary}
\newtheorem{conjecture}{Conjecture}
\newtheorem{pretheorem}[proposition]{Pretheorem}
\newtheorem{hypothesis}[proposition]{Hypothesis}
\newtheorem{example}[proposition]{Example}
\newtheorem{remark}[proposition]{Remark}
\newtheorem{ex}[proposition]{Exercise}
\newtheorem{cond}[proposition]{Conditions}
\newtheorem{cons}[proposition]{Construction}
\newcommand{\llabel}[1]{\label{#1}}
\newcommand{\comment}[1]{}
\newcommand{\sr}{\rightarrow}
\newcommand{\lr}{\longrightarrow}
\newcommand{\xr}{\xrightarrow}
\newcommand{\dw}{\downarrow}
\newcommand{\bdl}{\bar{\Delta}}
\newcommand{\zz}{{\bf Z\rm}}
\newcommand{\zq}{{\bf Z}_{qfh}}
\newcommand{\nn}{{\bf N\rm}}
\newcommand{\qq}{{\bf Q\rm}}
\newcommand{\nq}{{\bf N}_{qfh}}
\newcommand{\oo}{\otimes}
\newcommand{\uu}{\underline}
\newcommand{\ih}{\uu{Hom}}
\newcommand{\af}{{\bf A}^1}
\newcommand{\wt}{\widetilde}
\newcommand{\gm}{{\bf G}_m}
\newcommand{\dsr}{\stackrel{\sr}{\scriptstyle\sr}}
\newcommand{\PP}{$P_{\infty}$}
\newcommand{\tp}{\tilde{D}}
\newcommand{\HH}{$H_{\infty}$}
\newcommand{\ii}{\stackrel{\scriptstyle\sim}{\sr}}
\newcommand{\BB}{_{\bullet}}
\newcommand{\D}{\Delta}
\newcommand{\colim}{{\rm co}\hspace{-1mm}\lim}
\newcommand{\cf}{{\it cf} }
\newcommand{\msf}{\mathsf }
\newcommand{\mcal}{\mathcal }
\newcommand{\ep}{\epsilon}
\newcommand{\tl}{\widetilde}
\newcommand{\ub}{\mbox{\rotatebox{90}{$\in$}}}
\newcommand{\piece}{\vskip 3mm\noindent\refstepcounter{proposition}{\bf
\theproposition}\hspace{2mm}}
\newcommand{\subpiece}{\vskip 3mm\noindent\refstepcounter{equation}{\bf\theequation}\hspace{2mm}}{\vskip
3mm}
\numberwithin{equation}{subsection}
\begin{center}
{\Large\bf Cancellation theorem}\\
\vskip 4mm
{\large Vladimir Voevodsky}\\
{\em January 2002}
\end{center}
\vskip 4mm
\tableofcontents
\subsection{Introduction}
Let $SmCor(k)$ be the category of smooth finite correspondences over a
field $k$. Denote by $\gm$ the scheme $\af-\{0\}$. One defines the
sheaf with transfers $S^1_t$ by the condition that
$\zz_{tr}(\gm)=S^1_t\oplus\zz$ where $\zz$ is split off by the
projection to the point and the point $1$. For any scheme $X$ consider
the sheaf with transfers
$$F_X:Y\mapsto Hom(S^1_t\oo \zz_{tr}(Y), S^1_t\oo \zz_{tr}(X))$$
We prove that the obvious map $\zz_{tr}(X)\sr F_X$ defines a
quasi-isomorphism of singular simplicial complexes
$$C_*(\zz_{tr}(X))\sr C_*(F_X)$$
as complexes of presheaves i.e. for any $Y$ the map of complexes of
abelain groups
$$C_*(\zz_{tr}(X))(Y)\sr C_*(F_X)(Y)$$
is a quasi-isomorphism. We then deduce from this result the
cancellation theorem which asserts that if $k$ is a perfect field then
for any $K,L$ in $DM^{eff}_{-}(k)$ the map
$$Hom(K,K')\sr Hom(K(1),K'(1))$$
is bijective. 

This result was previously known in two particular situations. For
varieties over a field $k$ with resolution of singularities it was
proved in \cite{H3new}. For $K'$ being the motivic complex $\zz(n)[m]$
and any field $k$ it was proved in \cite{comparison}. Both proofs are
very long. 

The main part of our argument does not use the assumption that we work
with smooth schemes over a field and we give it for separated schemes
of finite type over a noetherian base. To be able to do it we define
in the first section the category of finite correspondences for
separated schemes of finite type over a base. The definition is a
straightforward generalization of the definition for schemes over a
field based on the constructions of \cite{SusVoe2new} and can be
skipped.  In the second section we define intersection of relative
cycles with Cartier divisors and prove the properties of this
construction which we need. In the third we prove our main theorem
\ref{main} and deduce from it the cancellation theorem over perfect
fileds \ref{main2}.

In this paper we say ``a relative cycle'' instead of ``an
equidimensional relative cycle''. All schemes are separated. The
letter $S$ is typically reserved for the base scheme which is assumed
to be noetherian. All the standard schemes ${\bf P}^1$, $\af$ etc. are
over $S$. When no confusion is possible we write $XY$ instead of
$X\times_S Y$.

\subsection{Finite correspondences}
For a scheme $X$ of finite type over a noetherian scheme $S$
we denote by $c(X/S)$ the group of finite relative cycles on $X$ over
$S$. In \cite{SusVoe2new} this group was denoted by
$c_{equi}(X/S,0)$. If $S$ is regular or if $S$ is normal and the
characteristic of $X$ is zero, $c(X/S)$ is the free abelian group
generated by closed irreducible subsets of $X$ which are finite over
$S$ and surjective over a connected component of $S$. For the general
definition see \cite[after Lemma 3.3.9]{SusVoe2new}. A morphism
$f:S'\sr S$ defines the pull-back homomorphism $c(X/S)\sr c(X
S'/S')$ which we denote by $cycl(f)$.

For two schemes $X,Y$ of finite type over $S$ we define the group
$c(X,Y)$ of {\em finite correspondences} from $X$ to $Y$ as
$c(X Y/X)$.  

Let us recall the following construction from \cite[\S
3.7]{SusVoe2new}. Let $X'\sr X\sr S$ be morphisms of finite type, ${\cal
W}$ a relative cycle on $X'$ over $X$ and ${\cal Z}$ a relative cycle
on $X$ over $S$. Then one defines a cycle $Cor({\cal W},{\cal Z})$ on
$X'$ as follows. Let $Z_i$ be the components of the support of $\cal
Z$ present with multiplicites $n_i$ and $e_i:Z_i\sr X$ the
corresponding closed embeddings. Let $e_i':Z_i\times_X X'\sr X'$
denote the projections. We set
$$Cor({\cal W},{\cal Z})=\sum_i n_i (e'_i)_*cycl(e_i)({\cal W})$$
where $(e'_i)_*$ is the (proper) push-forward on cycles.

Let $X$, $Y$ be schemes of finite type over $S$ and 
$$f\in c(X,Y)=c(X Y/X)$$
$$g\in c(Y,Z)=c(Y Z/Y)$$ 
finite correspondences. Let
$$p_X:X Y\sr Y$$
$$p_Y:X Y Z\sr X Z$$
be the projections. We define the composition $g\circ f$ by the
formula:
\begin{eq}\llabel{compdef}
g\circ f=(p_Y)_*Cor(cycl(p_X)(g),f)
\end{eq}
This operation is linear in both arguments and thus defines a
homomorphism of abelian groups
$$c(X,Y)\otimes c(Y,Z)\sr c(X,Z)$$
The lemma below follows immediately from the definition of $Cor(-,-)$
and the fact that the (proper) push-forward commutes with the
$cycl(-)$ homomorphisms (\cite[Prop. 3.6.2]{SusVoe2new}).
\begin{lemma}
\llabel{missing}
Let $Y\sr X\sr S$ be a sequence of morphisms of finite type,
$p:Y\sr Y'$ a morphism over $X$, ${\mathcal Y}\in
Cycl(Y/X,r)\oo\qq$ and ${\mathcal X}\in Cycl(X/S,s)\oo\qq$. Assume
that $p$ is proper on the support of ${\mathcal Y}$. Then 
$$p_*Cor({\mathcal
Y},{\mathcal X})=Cor(p_*({\mathcal Y},{\mathcal X})).$$
\end{lemma}
\begin{lemma}
\llabel{associativity}
For any $f\in c(X,Y)$, $g\in c(Y,Z)$, $h\in c(Z,T)$ one has
$$(h\circ g)\circ f=h\circ (g\circ f).$$
\end{lemma}
\begin{proof}
Consider the following diagram
$$
\begin{CD}
XT @<4<< XYT @>>> YT @. {} @. {}\\
@A7AA @A8AA @A2AA\\
XZT @<9<< XYZT @>>> YZT @>>> ZT @>>> T\\
@VVV @VVV @VVV @VVV @.\\
XZ @<5<< XYZ @>9>> YZ @>1>> Z @. {}\\
@. @VVV @VVV @. @.\\
{} @. XY @>3>> Y @. {} @. {}\\
@. @VVV @. @. @.\\
{} @. X
\end{CD}
$$
where the morphisms are the obvious projections. 
Note that all the squares are cartesian. We will also use the
projection $6:XZ\sr Z$.

We have $f\in c(XY/X)$, $g\in c(YZ/Y)$ and $h\in c(ZT/Z)$. The
compositions are given by:
$$g\circ f=5_*Cor(cycl(3)(g),f)$$
$$h\circ g=2_*Cor(cycl(1)(h),g)$$
$$(h\circ g)\circ f=4_*Cor(cycl(3)(h\circ
g),f)=4_*Cor(cycl(3)(2_*Cor(cycl(1)(h),g)),f)$$
$$h\circ (g\circ f)=7_*Cor(cycl(6)(h),g\circ f)=
7_*Cor(cycl(6)(h),5_*Cor(cycl(3)(g),f))$$
We have:
$$4_*Cor(cycl(3)(2_*Cor(cycl(1)(h),g)),f)=$$
$$=4_*Cor(8_*cycl(3)Cor(cycl(1)(h),g),f)=$$ 
$$=4_*8_*Cor(cycl(3)Cor(cycl(1)(h),g),f)=$$
$$=4_*8_*Cor(Cor(cycl(1\circ 9)(h),cycl(3)(g)),f)$$
where the first equality holds by \cite[Prop. 3.6.2]{SusVoe2new}, the
second by Lemma \ref{missing} and the third by
\cite[Th. 3.7.3]{SusVoe2new}.
We also have:
$$7_*Cor(cycl(6)(h),5_*Cor(cycl(3)(g),f))=$$
$$=7_*9_*Cor(cycl(6\circ 5)(h), Cor(cycl(3)(g),f))$$
by \cite[Lemma 3.7.1]{SusVoe2new}. We conclude that $(h\circ g)\circ
f=h\circ (g\circ f)$ by \cite[Prop. 3.7.7]{SusVoe2new}.
\end{proof}
We denote by $Cor(S)$ the category of finite correspondences whose
objects are schemes of finite type over $S$, morphisms are finite
correspondences and the composition of morphisms is defined by
(\ref{compdef}).

For a morphism of schemes $f:X\sr Y$ let $\Gamma_f$ be its graph
considered as an element of $c(XY/X)$. One verifies easily that
$\Gamma_{gf}=\Gamma_g \circ \Gamma_f$ and we get a functor $Sch/S\sr
Cor(S)$. Below we use the same symbol for a morphism of schemes and
its graph considered as a finite correspondence.

The external product of cycles defines pairings
$$c(X,Y)\oo c(X',Y')\sr c(X X', Y Y')$$
and one verifies easily using the results of \cite{SusVoe2new} that this
pairing extends to a tensor structure on $Cor(S)$ with $X\oo
Y:=X Y$.

\subsection{Intersecting relative cycles with divisors}
Let $X$ be a noetherian scheme and $D$ a Cartier divisor on $X$ i.e. a
global section of the sheaf ${\cal M}^*/{\cal O}^*$. One defines the
cycle $cycl(D)$ associated with $D$ as follows. Let $U_i$ be an open
covering of $X$ such that $D_{U_i}$ is of the form $f_{i,+}/f_{i,-}\in
{\cal M}^*(U_i)$. Then $cycl(D)$ is determined by the property that
$$cycl(D)_{|U_i}=cycl(f_{i,+}^{-1}(0))-cycl(f_{i,-}^{-1}(0))$$
where on the right hand side one considers the cycles associated with
closed subschemes (\cite[]{SusVoe2new}). One defines the support of $D$
as the closed subset $supp(D):=supp(cycl(D))$. 

We say that a cycle ${\cal Z}=\sum n_i z_i$ on $X$ intersects $D$
properly if the points $z_i$ do not belong to $supp(D)$. Let $Z_i$ be
the closure of $z_i$ considered as a reduced closed subscheme and
$e_i:Z_i\sr X$ the closed embedding. If ${\cal Z}$ and $D$ intersect
properly we define their intersection $({\cal Z},D)$ as the cycle
$$({\cal Z},D):=\sum n_i (e_i)_*(cycl(e_i^*(D)))$$
If $p:X\sr S$ is a morphism of finite type and ${\cal Z}$ is a
relative cycle of relative dimension $d$ over $S$, we say that $D$
intersects $\cal Z$ properly relative to $p$ (or properly over $S$) if
the dimension of fibers of $supp(D)\cap supp({\cal Z})$ over $S$ is
$\le d-1$. This clearly implies that ${\cal Z}$ intersects $D$
properly and $({\cal Z},D)$ is defined.
\begin{proposition}
\llabel{p1} Assume that $S$ is reduced and $D$ and $\cal Z$ intersect
properly over $S$. Then:
\begin{enumerate}
\item  $({\cal Z},D)$ is a relative cycle of relative
dimension $d-1$ over $S$
\item let $p:S'\sr S$ be a morphism and $p':X S'\sr X$ the
projection, assume that $(p')^*(D)$ is defined, then one has
\begin{eq}
\llabel{eqp1}
p^*(cycl(Z),D)=(p^*(cycl(Z)),(p')^*(D)).
\end{eq}
\end{enumerate}
\end{proposition}
\begin{proof}
One can easily see that $D$ intersects $\cal Z$ properly relative to
$p$ if and only if for any $i$ and any point $x$ of $X$ there exists a
neighborhood $U$ of $x$ such that $D_{|U}=f_+/f_{-}$ where $f_+$ and
$f_-$ are not zero divisors on the fibers of the maps $Z_i\sr S$.
Since the problem is local in $X$ we may assume that $D=D(f)$ where
$f$ is a regular function on $X$ which is not a zero divisor on fibers
of the maps $Z_i\sr S$. We write $({\cal Z},f)$ instead of $({\cal
Z},D(f))$.  

To prove both statements of the proposition it is clearly sufficient
to show that for any fat point
$Spec(k)\stackrel{x_0}{\sr}Spec(R)\stackrel{x_1}{\sr} S$ one has
$$(x_0,x_1)^*(({\cal Z},f))=((x_0,x_1)^*({\cal Z}), f  x'_1 
x'_0)$$
where $x_0', x_1'$ are the projections 
$$X Spec(k)\stackrel{x_0'}{\sr} X
Spec(R)\stackrel{x_1'}{\sr} X$$
Let $j:Spec(K(R))\sr Spec(R)$ be the generic point of $R$ and
$$pr:Spec(K(R)) X\sr X$$
the projection. Since $S$ is reduced, $Z_i$ are flat over a dense open
subset of $S$ and $(x_1j)^*(({\cal Z},f))=\sum m_j w_j$ is defined and
by Lemma \ref{flat} we have
$$(x_1j)^*(({\cal Z},f))=((x_1j)^*({\cal Z}), f\circ pr)$$
For a cycle ${\cal W}=\sum m_j w_j$ over the generic point of
$Spec(R)$ denote by $cl\,{\cal W}$ the same cycle considered over
$Spec(R)$. Since $cl(x_1j)^*({\cal Z})$ is a formal linear
combination of cycles of flat subschemes, Lemma \ref{flat} applies and
we conclude that:
$$(x_0,x_1)^*(({\cal Z},f))=x_0^*(cl(x_1j)^*(({\cal
Z},f)))=x_0^*(cl(x_1j)^*({\cal Z}), fx_1')=$$
$$=(x_0^*cl(x_1j)^*({\cal Z}),
fx_1'x_0')=((x_0,x_1)^*({\cal Z}),fx_1'x_0')$$
\end{proof}
\begin{lemma}
\llabel{two} Let $Z$ be a noetherian scheme and $f$ a regular function
on $Z$ which is not a zero divisor. Let further $e_i:Z_i\sr Z$ be the
connected components of $Z$ and $n_i$ their multiplicities. Then one
has:
\begin{eq}
\llabel{form1}
cycl(f^{-1}(0))=\sum n_i (e_i)_*(cycl(f_{|Z_i}^{-1}(0)))
\end{eq}
\end{lemma}
\begin{proof}
(P.Deligne) Replacing $Z$ by the local scheme of a generic point of
$f^{-1}(0)$ we may assume that $Z=Spec(A)$ is local of dimension $1$.  
By definition of the cycle associated to a closed subscheme
(\ref{form1}) is equivalent to 
$$l(A/(f))=\sum l(A/p_i) l(A/(p_i+(f)))$$
where $p_i$ are the minimal ideals of $A$ and $l(M)$ is the length of
an $A$-module. The function $M\sr l(M)$ extends to an additive
function on the $K_0$ of the abelain category of $A$-modules of finite
length. Consider the function $l_f:M\mapsto
l(cone(M\stackrel{f}{\sr}M))$. Equivalently we have
$l_f(M)=l(M\oo^{L}A/f)$ which shows that it is an additive function on
$K_0(A-mod)$. In $K_0(A-mod)$ one has
$$A=(\sum l(A/p_i) A/p_i)+M$$
where $M$ is a module of finite length. This implies:
$$l(A/f)=l_f(A)=\sum l(A/p_i)l_f(A/p_i)=\sum l(A/p_i)l(A/(p_i+(f))).$$
\end{proof}
\begin{lemma}
\llabel{flat} Let $Z$ be a closed subscheme of $X$ which is flat and
equidimensional over $S$ and $f$ be a regular function on $Z$ such
that $f^{-1}(0)\cap Z$ is equidimensional over $S$.  Then
$(cycl(Z),f)$ is a relative cycle and for any morphism $p:S'\sr S$ one
has:
$$p^*(cycl(Z),f)=(p^*(cycl(Z)),fp')$$
where $p':X S'\sr X$.
\end{lemma}
\begin{proof}
Lemma \ref{two} implies that
$$(cycl(Z),f)=cycl(f_{|Z}^{-1}(0)).$$ 
Our condition on $f$ implies that, for any $s\in S$, $f$ restricted to
the fiber $Z_s$ is not a zero divisor. Since $Z$ is flat over $S$ this
implies that $f_{|Z}^{-1}(0)$ is flat over $S$ (see
e.g. \cite[]{Miln}).  Since $Z$ is equidimensional of relative
dimension $d$, $f_{|Z}^{-1}(0)$ is equidimensional of relative
dimension $d-1$. We conclude by \cite[]{SusVoe2new} that
$cycl(f_{|Z}^{-1}(0))$ is a relative cycle.

For $p;S'\sr S$ we have:
$$p^*((cycl(Z),f))=p^*(cycl(f_{|Z}^{-1}(0)))=cycl((fp')_{|Z
S'}^{-1}(0))=$$
$$=(cycl(Z S'),fp')=(p^*(cycl(Z)),fp').$$
\end{proof}
\begin{cor}
\llabel{corr} Let $X'\stackrel{f}{\sr}X\sr S$ be morphisms of finite
type, ${\cal Z}$ a relative cycle on $X$ over $S$ and ${\cal W}$ a
relative cycle on $X'$ over $X$ of dimension $0$. Let further $D$ be a
Cartier divisor on $X'$ which intersects ${\cal W}$ properly over
$X$. Then $D$ intersects $Cor({\cal W},{\cal Z})$ properly over $S$
and one has:
\begin{eq}
\llabel{eqcorr}
(Cor({\cal W},{\cal Z}),D)=Cor(({\cal W},D),{\cal Z})
\end{eq}
\end{cor}
\begin{proof}
It is a straightforward corollary of the definition of
$Cor(-,-)$ and (\ref{eqp1}).
\end{proof}

\begin{lemma}
\llabel{pushfor} Let $f:X'\sr X$ be a morphism of schemes of finite
type over $S$, $\cal Z$ a relative cycle on $X'$ such that $f$ is
proper on $supp({\cal Z})$ and $D$ a Cartier divisor on $X$. Assume
that $f^*(D)$ is defined and $\cal Z$ intersects $f^*(D)$ properly
over $S$. Then $f_*({\cal Z})$ intersects $D$ properly over $S$ and
one has:
\begin{eq}
\llabel{eqp}
f_*({\cal Z},f^*(D))=(f_*({\cal Z}),D)
\end{eq}
\end{lemma}
\begin{proof}
Let $d$ be the relative dimension of $\cal Z$ over $S$. To see that
$f_*({\cal Z})$ intersects $D$ properly over $S$ we need to check that the
dimension of the fibers of $supp(D)\cap supp(f_*({\cal Z}))$ over $S$
is $\le d-1$. This follows from our assumption and the inclusion
$$supp(D)\cap supp(f_*({\cal Z}))\subset supp(D)\cap f(supp({\cal
Z}))=$$
$$=f(f^{-1}(supp(D))\cap supp({\cal Z}))=f(supp(f^*(D))\cap
supp({\cal Z}))$$
To verify (\ref{eqp}) it is sufficient to consider the situation
locally around the generic points of $f(supp(f^*(D))\cap supp({\cal
Z}))$. Therefore we may assume that $D=D(g)$ is the divisor of a
regular function $g$ and ${\cal Z}=z$ is just one point with the
closure $Z$. Replacing $X'$ by $Z$ and $X$ by $f(Z)$ we may assume
that $X$, $X'$ are integral, $f$ is surjective and $X$ is local of
dimension $1$. Let $A={\cal O}(X)$, $B={\cal O}(X')$. As in the proof
of \ref{two}, consider the function $l_g:M\mapsto l_A(M\oo^{L}A/g)$ on
$K_0(A-mod)$. This function vanishes on modules with the support in
the closed point which implies that 
$$l_g(B)=deg(f)l_g(A)=deg(f)l_A(A/g)$$ 
On the other hand $l_g(A)=l_A(B/(f^*(g)))$. Let $x'_i$ be the closed
points of $X'$, $k'_i$ their residue fields and $k$ the residue field
of the closed point of $X$. Let further $M_i$ be the part of
$B/(f^*(g))$ supported in $x_i'$. One can easily see that
$l_A(B/(f^*(g)))=\sum_i [k'_i:k]l_B(M_i)$. Combining our equalities we
get:
\begin{eq}
\llabel{cros}
deg(f)l_A(A/g)=\sum_i [k'_i:k]l_B(M_i)
\end{eq}
which is equivalent to (\ref{eqp}).
\end{proof}

\subsection{Cancellation theorem}
Consider a finite correspondence 
$${\cal Z}\in c(\gm X,\gm Y)=c(\gm X \gm Y/\gm X).$$
Let $f_1, f_2$ be the projections to the first and the second copy of
$\gm$ respectively and let $g_n$ denote the rational function
$(f_1^{n+1}-1)/(f_1^{n+1}-f_2)$ on $\gm X \gm Y$.
\begin{lemma}
\llabel{nex}
For any $\cal Z$ there exists $N$ such that for all $n\ge N$ the
divisor of $g_n$ intersects $\cal Z$ properly over $X$ and the cycle
$({\cal Z},D(g_n))$ is finite over $X$.
\end{lemma}
\begin{proof}
Let $\bar{f}_1\times \bar{q}:\bar{C}\sr {\bf P}^1 X$ be a finite
morphism which extends the projection $supp({\cal Z})\sr \gm
X$. Let $N$ be an integer such that the rational function
$\bar{f}^N_1/f_2$ is regular in a neighborhood of $\bar{f}_1^{-1}(0)$ and
the rational function $f_2/\bar{f}_1^{N}$ is regular in a neighborhood of
$\bar{f}_1^{-1}(\infty)$. Then for any $n\ge N$ one has:
\begin{enumerate}
\item the restriction of $g_nf_2$ to $supp({\cal Z})$ is regular on a
neighborhood of $\bar{f}_1^{-1}(0)$ and equals $1$ on $\bar{f}_1^{-1}(0)$
\item the restriction of $g_n$ to $supp({\cal Z})$ is regular a
neighborhood of $\bar{f}_1^{-1}(\infty)$ and equals $1$ on
$\bar{f}_1^{-1}(\infty)$
\end{enumerate}
Conditions (1),(2) imply that the divisor of $g_n$ intersects
$\cal Z$ properly over $X$ and that the relative cycle $({\cal
Z},D(g_n))$ is finite over $X$. 
\end{proof}
If $({\cal Z},D(g_n))$ is defined as a finite relative cycle we let
$\rho_n({\cal Z})\in c(X,Y)$ denote the projection of $({\cal
Z},D(g_n))$ to $X Y$.
\begin{remark}\llabel{4.2}\rm
Note that we can define a finite correspondence $\rho_g({\cal Z}):X\sr
Y$ for any function $g$ satisfying the conditions (1),(2) in the same
way as we defined $\rho_n=\rho_{g_n}$. In particular, if $n$ and $m$
are large enough then the function $tg_n+(1-t)g_{m}$ defines a finite
correspondence $h=h_{n,m}:X\af\sr Y$ such that
$h_{|X\times\{0\}}=\rho_m({\cal Z})$ and
$h_{|X\times\{1\}}=\rho_n({\cal Z})$, i.e. we get a canonical
$\af$-homotopy from $\rho_m({\cal Z})$ to $\rho_n({\cal Z})$.
\end{remark}
\begin{lemma}\llabel{str}
\begin{enumerate}[(i)]
\item For a finite correspondence ${\cal W}:X\sr Y$ and any $n\ge 1$
one has $\rho_n(Id_{\gm}\oo{\cal W})={\cal W}$
\item Let $e_X$ be the composition $\gm X\xr{pr} X\xr{\{1\}\times Id} \gm
X$. Then $\rho_n(e_X)=0$ for any $n\ge 0$.
\end{enumerate}
\end{lemma}
\begin{proof}
The cycle on $\gm X\gm Y$ over $\gm X$ which represents
$Id_{\gm}\oo{\cal W}$ is $\Delta_*(\gm\times {\cal W})$ where $\Delta$
is the diagonal embedding $\gm X Y\sr \gm X\gm Y$. The cycle
$(\Delta_*(\gm\times {\cal W}), g_n)$ is $\Delta_*(D\oo {\cal W})$
where $D$ is the divisor of the function $(t^{n+1}-1)/(t^{n+1}-t)$ on
$\gm$. The push-forward of $\Delta_*(D\oo {\cal W})$ to $XY$ is the cycle
$deg(D){\cal W}$. Since $deg(D)=1$ we get the first statement of the
lemma. 

The cycle $\cal Z$ on $\gm X\gm X$ representing $e_X$ is the image of the
embedding $\gm X\sr \gm X \gm X$ which is diagonal on $X$ and of the
form $t\mapsto (t,1)$ on $\gm$. This shows that the restriction of
$g_n$ to $supp({\cal Z})$ equals $1$ and $({\cal Z},D(g_n))=0$. 
\end{proof}
\begin{lemma}
\llabel{functor1} Let ${\cal Z}:\gm X\sr \gm Y$ be a finite
correspondence such that $\rho_n({\cal Z})$ is defined. Then for any
finite correspondence ${\cal W}:X'\sr X$, $\rho_n({\cal
Z}\circ({Id_{\gm}}\oo{\cal W}))$ is defined and one has
\begin{eq}
\llabel{eqf1}
\rho_n({\cal Z}\circ({Id_{\gm}}\oo{\cal W}))=\rho_n({\cal Z})\circ {\cal
W}
\end{eq}
\end{lemma}
\begin{proof}
Let us show that (\ref{eqf1}) holds. In the process it will become
clear that the left hand side is defined.  We can write $\rho_n({\cal
Z})\circ {\cal W}$ as the composition
$$X'\xr{\cal W} X\xr{({\cal Z},D(g_n))} \gm\gm Y\xr{pr} Y$$
and $\rho_n({\cal Z}\circ( {Id_{\gm}}\oo{\cal W}))$ as the composition
$$X'\xr{\cal Y}\gm\gm Y\xr{pr} Y$$
where ${\cal Y}=({\cal Z}\circ (Id_{\gm}\oo {\cal
W}),D(g_n))$. Consider the diagram
$$
\begin{CD}
\gm X'\gm Y @<p_1<< \gm X'X\gm Y @>>> \gm X\gm Y\\
@. @VVV @VVV\\
{} @. X'X @>p_2>> X\\
@. @VVV @.\\
{} @. X'
\end{CD}
$$
where the arrows are the obvious projections. If we consider $\cal Z$ as a
cycle of dimension $1$ over $X$ then the cycle ${\cal
Z}\circ({Id_{\gm}}\oo{\cal W})$, considered as a cycle over $X'$, is
$(p_1)_*Cor(cycl(p_2)({\cal Z}),{\cal W})$ and we have
$$((p_1)_*Cor(cycl(p_2)({\cal Z}),{\cal W}), D(g_n))=$$
$$=(p_1)_*(Cor(cycl(p_2)({\cal Z}),{\cal W}),
D(g_n))=(p_1)_*Cor((cycl(p_2)({\cal Z}), D(g_n)),{\cal W})=$$
$$=(p_1)_*Cor(cycl(p_2)({\cal Z},D(g_n)),{\cal
W})$$
where the first equality holds by (\ref{eqp}), the second by
(\ref{eqcorr}) and the third by (\ref{eqp1}). 

The last expression represents the composition ${\cal W}\circ ({\cal
Z},D(g_n))$ and we conclude that 
$$\rho_n({\cal Z})\circ {\cal W}=\rho_n({\cal Z}\circ( {Id_{\gm}}\oo{\cal
W}))$$
\end{proof}
\begin{lemma}
\llabel{functor3} Let ${\cal Z}:\gm X\sr \gm Y$ be a finite
correspondence such that $\rho_n({\cal Z})$ is defined. Then for any
morphism of schemes $f:X'\sr Y'$, $\rho_n({\cal Z}\oo f)$ is defined
and one has
\begin{eq}
\llabel{eqf3}
\rho_n({\cal Z}\oo f)=\rho_n({\cal Z})\oo f
\end{eq}
\end{lemma}
\begin{proof}
Consider the diagram
$$
\begin{CD}
\gm X X'\gm Y Y' @<p_1<< \gm X X'\gm Y @>>> \gm X \gm Y\\
@. @VVV @VVV\\
{} @. X X' @>p_2>> X
\end{CD}
$$
where $p_1$ is defined by the embedding $X'\xr{f\times Id} X'Y'$ and
the rest of the morphisms are the obvious projections. Consider $\cal
Z$ as a cycle over $X$. Then $\rho_n({\cal Z}\oo f)$ is given by the
composition
$$\gm X X' \xr{{\cal Y}_1} \gm\gm Y\xr{pr} YY'$$
where ${\cal Y}_1=((p_1)_*cycl(p_2)({\cal Z}),g_n)$ and $\rho_n({\cal
Z})\oo f$ by the composition
$$\gm X X' \xr{{\cal Y}_2} \gm\gm Y\xr{pr} YY'$$
where ${\cal Y}_2=(p_1)_*(cycl(p_2)(({\cal Z},g_n)))$. The equality
${\cal Y}_1={\cal Y}_2$ follows from (\ref{eqp}) and (\ref{eqp1}).
\end{proof} 
For our next result we need to use presheaves with transfers. A
presheaf with transfers on $Sch/S$ is an additive contravariant
functor from $Cor(S)$ to the category of abelian groups.  For $X$ in
$Sch/S$ we let $\zz_{tr}(X)$ denote the functor represented by $X$ on
$Cor(S)$. One defines tensor product of presheaves with transfers in
the usual way such that $\zz_{tr}(X)\oo\zz_{tr}(Y)=\zz_{tr}(X\times
Y)$. To simplify notations we will write $X$ instead of $\zz_{tr}(X)$
and identify morphisms $\zz_{tr}(X)\sr\zz_{tr}(Y)$ with finite
correspondences $X\sr Y$. Note in particular that $\gm$ denotes the
presheaf with transfers $\zz_{tr}(\gm)$ not the presheaf with
transfers represented by $\gm$ as a scheme. To preserve compatibility
with the notation $XY$ for the product of $X$ and $Y$ we write $FG$
for the tensor product of presheaves with transfers $F$ and $G$.

Let $S^1_t$ denote the presheaf with transfers $ker(\gm\sr
S)$. We consider it as a direct summand of $\gm$ with
respect to the projection $Id-e$ where $e$ is defined by the
composition $\gm\sr S\xr{1}\gm$. In the following theorem we let
$f\cong g$ denote that the morphisms $f$ and $g$ are
$\af$-homotopic. 
\begin{theorem}
\llabel{main} Let $F$ be a presheaf with transfers such that there is
an epimorphism $X\sr F$ for a scheme $X$. Let $\phi:S^1_t\oo F\sr
S^1_t  Y$ be a morphism. Then there exists a unique up to an
$\af$-homotopy morphism $\rho(\phi):F\sr Y$ such that $Id_{S^1_t}\oo
\rho(\phi)\cong \phi$.
\end{theorem}
\begin{proof}
Let us fix an epimorphism $p:X\sr F$. Then the morphism $\phi$ defines
a finite correspondence ${\cal Z}:\gm X\sr \gm Y$ and for $n$
sufficiently large we may consider $\rho_n({\cal Z}):X\sr Y$. Lemma
\ref{functor1} implies immediately that $\rho_n({\cal Z})$ vanishes on
$ker(p)$ and therefore it defines a morphism $\rho_n(\phi):F\sr X$.

Consider a morphism $\phi$ of the form $Id_{S^1_t}\oo \psi$. Then
$\cal Z$ is of the form $(Id_{\gm}-e)\oo {\cal W}$ where ${\cal
W}:X\sr Y$ corresponds to $\psi$. By Lemma \ref{str} we have
$\rho_n({\cal Z})={\cal W}$ and therefore
$\rho_n(Id_{S^1_t}\oo\psi)=\psi$ for any $n\ge 1$. If $\rho, \rho'$
are two morphims such that $Id_{S^1_t}\oo\rho\cong \phi$ and
$Id_{S^1_t}\oo\rho'\cong \phi$ then for a sufficiently large $n$ we have
$$\rho=\rho_n(Id_{S^1_t}\oo\rho)\cong \rho_n(Id_{S^1_t}\oo\rho')=\rho'$$
This implies the uniqueness part of the theorem. 

To prove the existence let us show that for a sufficiently large $n$
one has $Id_{S^1_t}\oo \rho_n(\phi)\cong \phi$. Let $\wt{\phi}$ be the
morphism $\gm F\sr \gm Y$ defined by $\phi$ and let
$$\wt{\phi}^*:F  \gm\sr Y  \gm$$
be the morphism obtained from $\wt{\phi}$ by the obvious
permutation. 
\begin{lemma}
\llabel{perm} The morphisms $\wt{\phi}\oo(Id_{\gm}-e)$ and
$(Id_{\gm}-e)\oo\wt{\phi^*}$ are $\af$-homotopic.
\end{lemma}
\begin{proof}
One can easily see that these two morphisms are obtained from the
morphisms
$$\phi\oo Id_{S^1_t}, Id_{S^1_t}\oo \phi^*:S^1_t  F  S^1_t\sr S^1_t 
Y  S^1_t$$
by using the standard direct sum decomposition. One can see further
that $\phi\oo Id_{S^1_t}=\sigma_Y (Id_{S^1_t}\oo \phi^*) \sigma_F$ where
$\sigma_F$ and $\sigma_Y$ are the permutations of the two copies of
$S^1_t$ in $S^1_t  F  S^1_t$ and $S^1_t  Y  S^1_t$
respectively. Lemma \ref{expl} below implies now that $\phi\oo
Id_{S^1_t}\cong Id_{S^1_t}\oo\phi^*$.
\end{proof}
\begin{lemma}
\llabel{expl}
The permutation on $S^1_t  S^1_t$ is $\af$-homotopic to $\{-1\}Id\oo
Id$ where $\{-1\}:S^1_t\sr S^1_t$ is defined by the morphism
$\gm\xr{x\mapsto x^{-1}}\gm$. 
\end{lemma}
\begin{proof}
The same arguments as the ones used in \cite[p.142]{SusVoe3} show that
for any scheme $X$ and any pair of invertible functions $f,g$ on
$X$ the morphism $X\xr{f\oo g}S^1_t  S^1_t$ is $\af$-homotopic to
the morphism $g\oo f^{-1}$. This implies immediately that the permutation
on $S^1_t  S^1_t$ is $\af$-homotopic to the morphism
$Id\oo(\{-1\}Id)$ where $\{-1\}Id:S^1_t\sr S^1_t$ is the morphism
defined by the map $\gm\xr{x\mapsto x^{-1}}\gm$. 
\end{proof}
For a sufficiently large $n$ we have
$$\rho_n(\phi\oo(Id_{\gm}-e))=\rho_n(\phi)\oo
(Id_{\gm}-e)$$
by Lemma \ref{functor3}. On the other hand
$$\rho_n((Id_{\gm}-e)\oo\phi^*)=\phi^*$$
by Lemma \ref{str}. By Lemma \ref{perm} we conclude that 
$$\phi^*\cong \rho_n(\phi)\oo
(Id_{\gm}-e)$$
which is equivalent to  $Id_{S^1_t}\oo \rho_n(\phi)\cong \phi$. Theorem
\ref{main} is proved. 
\end{proof}
\begin{cor}
\llabel{main1}
Denote by $F_Y$ the presheaf 
$$X\mapsto
Hom(S^1_t X,S^1_t Y)$$
and consider the obvious map $Y\sr F_Y$. Then for any $X$
the corresponding map of complexes of abelian groups
$$C_*(Y)(X)\sr C_*(F_Y)(X)$$
is a quasi-isomorphism 
\end{cor}
\begin{proof}
Let $\Delta^n\cong {\bf A}^n$ be the standard algebraic simplex and
$\partial\Delta^n$ the subpresheaf in $\Delta^n$ which is the
union of the images of the face maps $\Delta^{n-1}\sr
\Delta^n$. Then the n-th homology group of the complex
$C_*(F)(X)$ for any $F$ is the group of homotopy classes of maps from
$X \oo (\Delta^n/\partial\Delta^n)$ to $F$. Our result
now follows directly from \ref{main}.
\end{proof}
\begin{cor}
\llabel{main2} Let $k$ be a perfect field. Then for any $K,L$ in
$DM^{eff}_{-}(k)$ the map $Hom(K,L)\sr Hom(K(1),L(1))$ is a bijection.
\end{cor}
\begin{proof}
Since $DM^{eff}_{-}$ is generated by objects of the form $X$
it is enough to check that for smooth schemes $X,Y$ over $k$ and
$n\in\zz$ one has
$$Hom(S^1_t  X, S^1_t 
Y[n])=Hom(X,Y[n])$$
By Corollary \ref{main1} we know that the map 
$$Y\sr F_Y=\uu{Hom}(S^1_t, S^1_t  Y)$$
is an isomorphism in $DM$. Let us show now that for any sheaf with
transfers $F$ and any $X$ one has
\begin{eq}
\llabel{eqal1}
Hom_{DM}(S^1_t  X, F[n])=Hom_{DM}(X,
\uu{Hom}(S^1_t, F)[n])
\end{eq}
The left hand side of (\ref{eqal1}) is the hypercohomology group ${\bf
H}^n(\gm X, C_*(F))$ modulo the subgroup ${\bf H}^n(X,
C_*(F))$. The right hand side is the hypercohomology group ${\bf
H}^n(X, C_*\uu{Hom}(\gm , F))$ modulo similar
subgroup. Let $p:\gm X\sr X$ be the projection. It
is easy to see that (\ref{eqal1}) asserts that ${\bf
R}p_*(C_*(F))\cong C_*(p_*(F))$. There is a spectral sequence which
converges to the cohomology sheaves of ${\bf R}p_*(C_*(F))$ and starts
with the higher direct images $R^ip_*(\uu{H}^j(C_*(F)))$. We need to
verify that $R^ip_*(\uu{H}^j(C_*(F)))=0$ for $i>0$ and that
$p_*(\uu{H}^j(C_*(F)))=\uu{H}^j(C_*(p_*(F)))$. Both statements follow
from \cite[Prop. 4.34, p.124]{H2new} and the comparison of Zariski and
Nisnevich cohomology for homotopy invariant presheaves with
transfers.
\end{proof}

\comment{
\subsection{Concluding remarks}
\piece\llabel{placeofdm} Consider the category
$\Delta^{op}Shv_{Nis}Cor(S)$ of simplicial Nisnevich sheaves with
transfers on $Sch/S$. It is an analog of the category of pointed
sheaves of sets on $Sch/S$ and one can use this analogy to define the
{\em unstable $\af$-homotopy category of finite 
correspondences} $H_{\af}(Cor(S))$. This is the localization of
$\Delta^{op}Shv_{Nis}Cor(S)$ by weak equivalences which can be defined
either through $\af$-local objects or through $\bdl$-closures. One can
also think of $H_{\af}(Cor(S))$ in terms of complexes using the
equivalence 
$$\Delta^{op}Shv_{Nis}Cor(S)\cong Compl_{\ge 0}Shv_{Nis}Cor(S).$$ 
\piece\llabel{susp}
As for schemes we have two circles $S^1_s$ and $S^1_t$ where
$S^1_s=cone(\zz\sr 0)$ is the simplicial circle and
$S^1_t=\zz_{tr}(\gm)/\zz_{tr}(\{1\})$ is the circle discussed
above. The tensor product with the circles gives the suspension
functors $\Sigma^1_s$ and $\Sigma^1_t$. The s-suspension coincides
with the simplicial suspension or, in the language of complexes, with
the shift $K\mapsto K[1]$ and preserves $\af$-equivalences. The 
t-suspension preserves $\af$-equivalences between objects which are
termwise direct sums of sheaves of the form $\zz_{tr}(X)$. Using
appropriate resolution functor $G$ one can define ${\bf L}\Sigma^1_t$
as $\Sigma^1_tG$. 

\piece\llabel{stable} We can use s- and t-suspensions to define the
corresponding stable homotopy categories $SH_s(Cor(S))$,
$SH_t(Cor(S))$ and $SH_T(Cor(S))$ where the last one is obtained by
stablization in both s- and t-. In this language, the category $D(S)$
considered above is the s-stable homotopy category of {\em connective
s-spectra}. If we were to consider complexes unbounded in both
directions in the definition of $D(S)$ we would get the whole s-stable
category.

\piece\llabel{susp2} For smooth schemes over a perfect field the
unstable homotopy category $H_{\af}(SmCor(k))$ is almost stable with
respect to both s- and t-supensions: $\Sigma^1_s$ and $\Sigma^1_t$ are
full embeddings. The s-suspension part is a direct corollary of the
main theorem about homotopy invariant sheaves with transfers
(\cite[]{H2new}). Modulo the s-suspension part the t-suspenison part is
our \ref{main2}.  Below we give two examples which show that
over more general schemes this is no longer true. We use an old
example from \cite{} to show that the s-suspension is not a full
embedding over {\em any} scheme of positive dimension. We further give
an example of a singular local henselian scheme $S$ of dimension $2$
such that the t-suspension is not a full embedding on the s-stable
homotopy category of finite correspondences over $S$. This clearly
implies that it is not a full embeding on the unstable category. I
expect that the t-suspension is a full embedding on the s-stable
category over any {\em regular} scheme but I do not know what to
expect from the t-suspension on the unstable category.
\begin{example}\rm
Let $S$ be any scheme $j:U\subset S$, $i:Z\sr S$ be a complimentary
pair of open and closed embeddings such that $U,Z\ne \emptyset$. The
functor $j^*:Shv_{Nis}(Cor(S))\sr Shv_{Nis}(Cor(U))$ has a left
adjoint $j_!$ which takes $\zz_{tr}(X\sr U)$ to $\zz_{tr}(X\sr S)$. The
functor $i^*$ has a right adjoint $i_*$. The composition of two
adjunctions $j_!j^*\sr Id\sr i_*i^*$ is zero which gives us a morphism
$j_!j^*\sr ker(Id\sr i_*i^*)$. Consider
the sheaf ${\cal O}^*$ as a sheaf with transfers. Then the morphism
$$j_!j^!({\cal O}^*)\sr ker ({\cal O}^*\sr i_*i^*({\cal O}^*))$$ is an
s-stable $\af$-weak equivalence but not an $\af$-weak equivalence. To
see the later it is sufficient to notice that both sheaves are
$\af$-local and that this morphism is not an isomorphism. To see the
former one can use the reasoning similar to the one used in the proof
of the glueing theorem \cite[Th. 3.2.21]{MoVo} to show that the
cokernel of this morphism is $\af$-weakly equivalent to zero.
\end{example}
\piece\llabel{7.1.5}
(short ``exact'' sequence for $G_m$)
\piece\llabel{7.2} For a morphism $p:X\sr S$ of finite type we let
$$p_{\#}:Shv(Cor(X))\sr Shv(Cor(S))$$ 
denote the functor left adjoint to $p^*$. It takes $\zz_{tr}(Y\sr X)$
to $\zz_{tr}(Y\sr S)$. If $p$ is an open embedding $p_{\#}$ coincides
with the functor usually denoted by $p_!$. Let $Z\sr S$ be a closed
embedding and $p:X\sr S$ be an abstract blow-up with center in $Z$
i.e. a proper morphism such that $p^{-1}(S-Z)\sr S-Z$ is an
isomorphism. Let further $X\stackrel{i}{\sr}Y\stackrel{q}{\sr} S$ be
any decomposition of $p$ into a closed embedding and a smooth
morphism. Let $j^X:Y-i(X)\sr Y$, $j^Z:Y-i(p^{-1}(Z))\sr Y$ be open
embeddings. Using \ref{7.1.5}, one concludes easily that in the
s-stable homotopy category of finite correspondences one has:
$$Hom_S(X/p^{-1}(Z),\Sigma^m_s{\bf G}_m)=Hom_{Y}(\zz,
\Sigma^m_sj^Z_!({\bf G}_m)/j^X_!({\bf G}_m))$$
Applying $q_{\#}$ we get a map from this group to 
$$Hom_S(Y,\Sigma^m_s\Sigma^1_t(Y-i(p^{-1}(Z)))/(Y-i(X)))$$  
Under our assumption $X-p^{-1}(Z)\cong S-Z$ in particular it is smooth
over $S$ and if $n=dim(Y/S)$ we get by the purity theorem that 
$$(Y-i(p^{-1}(Z)))/(Y-i(X))=(Y-i(p^{-1}(Z)))/((Y-i(p^{-1}(Z)))-i(X-p^{-1}(Z)))\cong$$
$$\cong \Sigma^n_T(S-Z)$$
Hence, we get a map:
$$Hom_S(X/p^{-1}(Z),\Sigma^m_s{\bf G}_m)\sr
Hom_S(Y,\Sigma^m_s\Sigma^1_t\Sigma^n_T(S-Z))$$
If $Z=\emptyset$ and $X=S$ this map is just the composition with the
Gysin map $Y\sr Y/(Y-X)\cong \Sigma^n_T(X)$. 

For a closed subset $W\subset X$ denote the group of morphisms in the
s-stable homotopy category of finite correspondences
$Hom(X/W,\Sigma^1_s\Sigma^1_t)$ by $Pic^h(X,W)$. There is a map from
the usual relative Picard group to $Pic^h$ and if both $X$ and $W$ are
smooth this map is an isomorphism. Applying our construction for
$Y={\bf P}^n_S$ and $m=1$ we get a map
$$Pic^h(X,p^{-1}(Z))\sr Hom({\bf P}^n, \Sigma_T^{n+1}(S-Z))$$
Composing with the fundamental class $T^n\sr {\bf P}^n$ and the
inclusion $S-Z\sr S$ we finally get:
$$Pic^h(X,p^{-1}(Z))\sr Hom_S(T^n, T^{n+1})$$
One verifies easily that this map commutes with pull-backs and that
for $Z=\emptyset$ the resulting map 
$$Hom(S,T)\sr Hom_S(T^n, T^{n+1})$$
is just the n-fold suspension.

For a morphism of schemes $f:C\sr \gm $
let $\tilde{f}$ be the corresponding morphism $\zz_{tr}(C)\sr S^1_t$
in $D(S)$. Then one has $\widetilde{ff'}=\tilde{f}+\tilde{f}'$ where
$+$ sign refers to the sum of morphisms in the additive category
$D(S)$ and the product $ff'$ is the product of invertible functions.

$$
\begin{CD}
A @>>> B\\
@| @|
\end{CD}
$$

$$A\xrightarrow[aaaaaaaaaa]{b}X$$
}
\end{document}